\documentclass[11pt]{article}
\usepackage{pictex}
\usepackage{amssymb}

\renewcommand{\baselinestretch}{1.2}
\newtheorem{thm}{Theorem}
\newtheorem{cla}{Claim}
\newtheorem{lem}{Lemma}
\newtheorem{rem}{Remark}
\newtheorem{pro}{Proposition}

\newtheorem{cor}{Corollary}

\newcommand\proof{\noindent {\bf Proof:\hspace{.05in}}}
\newcommand\qed{\hfill $\Box$ \smallskip}

\begin{document}

\title{\bf Maximum independent sets near the upper bound}
\renewcommand\baselinestretch{1}\small

\author{
{\sl Ingo Schiermeyer}\\
\small Institut f\"ur Diskrete Mathematik und Algebra\\
\small Technische Universit\"at Bergakademie Freiberg
\\ \small 09596 Freiberg,
Germany\\
\small Ingo.Schiermeyer@tu-freiberg.de}
%{\sl Arnfried Kemnitz} \\ Computational Mathematics\\
%Technische Universit\"at
%Braunschweig\\ 38023 Braunschweig, Germany\\
%a.kemnitz@tu-bs.de\\
%\and
%{\sl Gabriel Semani\v{s}in}\\
%Mirko Hor\v{n}\'ak, , Roman Sot\'ak, \{Mirko.Hornak, , Roman.Sotak\}
%\thanks{This work was supported by the Slovak Science and Technology
%Assistance Agency under the contract APVV-0023-10 and Slovak VEGA
%grant 1/0652/12.
%}
%\small Institute of Mathematics\\
%\small P.J. \v{S}af\'arik University in Ko\v{s}ice
%\\ \small 04001 Ko\v{s}ice, Slovakia\\
%\small Gabriel.Semanisin@upjs.sk}

\date{\today}
\maketitle

\begin{abstract}

The size of a largest independent set of vertices in a given graph $G$ is denoted by $\alpha(G)$ and is called its {\it independence number} (or {\it stability number}).
Given a graph $G$ and an integer $K,$ it is NP-complete to decide whether $\alpha(G) \geq K.$ An upper bound for the independence number $\alpha(G)$
of a given graph $G$ with $n$ vertices and $m $ edges is given by \\$\alpha(G) \leq p:=\lfloor\frac{1}{2} + \sqrt{\frac{1}{4} + n^2 - n - 2m}\rfloor$.

In this paper we will consider maximum independent sets near this upper bound. Our main result is the following:
There exists an algorithm with time complexity $O(n^2)$ that, given as an input a graph $G$ with $n$ vertices, $m$ edges,
$p:=\lfloor\frac{1}{2} + \sqrt{\frac{1}{4} + n^2 - n - 2m}\rfloor,$ and an integer $k \geq 0$ with $p \geq 2k+1,$ returns an induced subgraph $G_{p,k}$ of $G$ with $n_0 \leq p+2k+1$ vertices such that
$\alpha(G) \leq p-k$ if and only if $\alpha(G_{p,k}) \leq p-k.$ Furthermore, we will show that we can decide in time $O(1.2738^{3k} + kn)$ whether $\alpha(G_{p,k}) \leq p-k.$

\end{abstract}

\renewcommand\baselinestretch{1.2}\normalsize

\section{Introduction}

We use \cite{W} for terminology and notation not defined here and
consider finite and simple graphs only. Given  two graphs $G$ and $H,$
then $G \cup H$ and $G + H$ denote the union and the join of $G$ and $H,$ respectively.

An \emph{independent set} (also called \emph{stable set}) in a graph $G$ is a subset of pairwise
non-adjacent vertices. The size of a largest independent set of vertices in a given graph $G$ is denoted by $\alpha(G)$ and is called its {\it independence number} (or {\it stability number}).
The {\sc Maximum Independent Set} problem (MIS for short) asks for finding
an independent set of maximum cardinality in a graph. This problem is known to be NP-hard in general.
Moreover, it remains NP-hard under several restrictions, for example for triangle-free graphs
\cite{Po}, $K_{1,4}$-free graphs \cite{Mi}, planar graphs of maximum degree at most three
\cite{GaJo}, and for sparse or dense graphs \cite{S}. On the other hand, for some subclasses of these graph classes
the MIS problem is solvable in polynomial time, for example
for graphs with maximum degree two (folklore), for subclasses of subcubic graphs \cite{BLS}, and for $K_{1,3}$-free graphs \cite{Mi,Sb}.

%\medskip

%For each fixed pair $\alpha, c > 0$ let INDEPENDENT SET $(m \leq cn^{\alpha})$ and INDEPENDENT SET $(m \geq {n \choose 2} -  %cn^{\alpha})$ be the problem
%INDEPENDENT SET restricted to graphs with $m \leq cn^{\alpha}$ or $m \geq {n \choose 2} -  cn^{\alpha}$ edges, respectively.

%\begin{thm}
%INDEPENDENT SET $(m \leq cn^{\alpha})$ is NP-complete.
%\end{thm}

%Note that INDEPENDENT SET $(m \leq k)$ can be solved in time $O(n^k)$ for any fixed integer $k \geq 1.$

%\begin{thm}
%INDEPENDENT SET $(m \geq {n \choose 2} -  cn^{\alpha})$ is NP-complete.
%\end{thm}

%Note that INDEPENDENT SET $(m \geq {n \choose 2} -  k)$ can be solved in time $O(n^k)$ for any fixed integer $k \geq 1.$

\section{Algorithmic bounds for the independence number of a graph}

Various lower and upper bounds for the independence number of a graph in terms of other graph parameters have
been shown in the last decades. A survey on these bounds is given in \cite{L}, where the number of lower bounds is considerably
larger than the number of upper bounds.

Recently, Dvo\v{r}\'ak and Lidick\'y have considered maximum independent sets near the {\it lower bound} in a graph,
where the lower bound is $\frac{n}{\Delta}.$

\medskip

As mentioned in \cite{DL}, by Brook's Theorem \cite{B}, every graph of order $n$, maximum degree at most $\Delta \geq 3$
and clique number at most $\Delta$ is $\Delta$-colourable, and thus has an independent set of size at least $\frac{n}{\Delta}.$

\begin{thm}(Dvo\v{r}\'ak and Lidick\'y)\label{t1}
There exists an algorithm with time complexity $O(\Delta^2 n)$ that, given as an input an integer
$\Delta \geq 3,$ a graph $G$ of order $n$ with $max(\Delta(G), \omega(G)) \leq \Delta,$ and an integer $k \geq 0,$ returns an induced subgraph $G_0$
of $G$ with $n_0 \leq 114\Delta^3 k$ vertices such that $\alpha(G) \geq \frac{n}{\Delta} +k$ if and only if $\alpha(G_0) \geq \frac{n_0}{\Delta} +k.$
\end{thm}

Hence the problem of deciding whether such a graph has an independent set of size at least $\frac{n}{\Delta} +k$ has a kernel of size $O(k).$ Note that these graphs belong to the class of sparse graphs. As mentioned in \cite{DL}, such an instance can be solved by brute force leading to an $2^{O(\Delta^3k)} + O(\Delta^2n)$ algorithm.

\medskip

In this paper we will consider independent sets near the upper bound in general graphs, where the $\it upper \ bound$ is the following bound (cf. \cite{L}).

\begin{thm}\label{t2}
Let $G$ be a graph of order $n$ and size $m.$ Then
$$\alpha(G) \leq p:=\lfloor\frac{1}{2}+\sqrt{\frac{1}{4}+n^2-n-2m}\rfloor.$$
\end{thm}

%By the definition, $m(\overline{G}) < {p+1 \choose 2}.$

The main result in this paper is the following analogue of Theorem \ref{t1}.

\begin{thm}\label{t3}
There exists an algorithm with time complexity $O(n^2)$ that, given as an input a graph $G$ of order $n$ with $p:=\lfloor\frac{1}{2}+\sqrt{\frac{1}{4}+n^2-n-2m}\rfloor,$ and an integer $k \geq 0$ with $p \geq 2k+1,$ returns an induced subgraph $G_{p,k}$
of $G$ with $n_0 \leq p+2k+1$ vertices such that $\alpha(G) \leq p - k$ if and only if $\alpha(G_{p,k}) \leq p - k.$
\end{thm}

Furthermore, we will show that we can decide in time $O(1.2738^{3k} + kn)$ whether $\alpha(G_{p,k}) \leq p-k.$

\section{Auxiliary results}

Given a graph $G$ of order $n,$ a subset $U \subset V(G)$ of vertices is called a {\it vertex cover} of $G,$ if for each edge of $G$ at least one of its incident vertices belongs to $U.$ The following property is well-known:

\centerline{$U$ is a vertex cover of $G$ if and only if $V(G) \setminus U$ is an independent set of $G$}

\medskip

Hence the problem of computing a minimum vertex cover is also NP-hard. However, given an integer $k \geq 1,$ the decision problem
whether $G$ has a vertex cover of size at most $k$ is {\it fixed parameter tractable}. This means, there exists an algorithm with time complexity $f(k) \cdot p(n),$ where $f(k)$ is a computable function depending only on $k$ and $p(n)$ is a polynomial depending on $n.$ Several fpt-algorithms have been developed for the vertex cover problem, where the currently fastest algorithm is due to
Chen et al. \cite{CKX} with time complexity $O(1.2738^k + kn).$

%Hence the problem of deciding whether such a graph has an independent set of size at most $p-k$ has a kernel of size $O(k).$
%Note that these graphs belong to the class of dense graphs.

%Such an instance can be solved by brute force in time $O(2^{2pk}) + O(n).$ Alternatively, we can ask whether $G_0$ has a vertex %cover of size at most $t = n_0 - (p-k+1).$
%Two prominent fpt-algorithms for solving VERTEX COVER with input $t$ are the algorithms of Buss with time complexity
%$O(tn + 2^t t^{2t+2})$ and of Papadimitriou and Yannakakis with time complexity $O(3^ttn).$
%For simplyfication, let $t:= (2k-1)(p+1) > (2k-1)p+k-1 = 2pk - (p-k+1) > n_0. $ Then these two algorithms
%could be applied with time complexity $O((2k-1)(p+1)n_0 + 2^{(2k-1)(p+1)})((2k-1)(p+1))^{2(2k-1)(p+1)}$ and $O(3^{(2k-1)(p+1)(2k-1)(p+1)n_0}),$ respectively.

\section{Main results}

For a given a graph $G$ with $n$ vertices, $m$ edges, $p:=\lfloor\frac{1}{2} + \sqrt{\frac{1}{4} + n^2 - n - 2m}\rfloor,$ and an integer $k \geq 0$ with $p \geq 2k+1,$
let $H \subset G$ be the subgraph of $G$ induced by all vertices having degree at least $n-p+k.$ Let $G_{p,k}=G - H.$

\begin{thm}\label{t4}
There exists an algorithm with time complexity $O(n^2)$ that, given as an input a graph $G$ of order $n$ with $n$ vertices, $m$ edges, $p:=\lfloor\frac{1}{2}+\sqrt{\frac{1}{4}+n^2-n-2m}\rfloor,$ and an integer $k \geq 0$ with $p \geq 2k+1,$ returns an induced subgraph $G_{p,k}$
of $G$ with $n_0 \leq p+2k+1$ vertices such that $\alpha(G) \leq p - k$ if and only if $\alpha(G_{p,k}) \leq p - k.$
\end{thm}

\proof
If $\alpha(G) \leq p - k$ then $\alpha(G_{p,k}) \leq p - k$ since $G_{p,k}$ is an induced subgraph of $G.$ Conversely, suppose $\alpha(G_{p,k}) \leq p - k,$
but $\alpha(G) \geq p-k+1.$ Let $I$ be an independent set with $p-k+1$ vertices and $H \subset G$ be the subgraph of $G$ induced by all vertices
having degree at least $n-p+k.$ Then $I \cap V(H) = \emptyset.$ Hence $I \subset G_{p,k} = G - H.$ Then $d_G(v) \leq n-(p-k+1) = n-p+k-1$ for all vertices in $G_{p,k}.$ Using the handshaking lemma
we obtain

$$\frac{n_0}{2}(p-k) = \frac{n_0}{2}(n-1-(n-p+k-1)) \leq \frac{1}{2}\sum_{v \in V(G_0)}(n-1-d_G(v)) = m(\overline{G}) < {p+1 \choose 2}$$

\noindent implying $n_0 < \frac{p(p+1)}{p-k}.$ Using $\frac{p(p+1)}{p-k} = p + \frac{p(k+1)}{p-k}$ we obtain $n_0 \leq p+2k+1$ for all $p \geq 2k.$
\qed

\medskip

Now we are searching for a vertex cover of size $t(k) = n_0 - (p-k+1).$ For all $p \geq 2k$ we obtain
$$t(k) < \frac{p(p+1)}{p-k} - (p-k+1) = \frac{2kp-k^2+k}{p-k} = 2k + \frac{k(k+1)}{p-k} \leq 3k+1.$$

\noindent Hence $t(k) \leq 3k.$

\medskip

Using the vertex cover fpt-algorithm of Chen et al. we obtain the next theorem.

\begin{thm} There exists an $O(1.2738^{3k} + kn)$ fpt-algorithm (fixed parameter tractable) to decide whether $\alpha(G_{p,k}) \leq p-k.$
\end{thm}

\begin{cor} Given as an input a graph $G$ with $n$ vertices, $m$ edges,
$p:=\lfloor\frac{1}{2} + \sqrt{\frac{1}{4} + n^2 - n - 2m}\rfloor,$ and an integer $k \geq 0$ with $p \geq 2k+1,$
it can be decided in time $O(1.2738^{3k} + kn)$ whether $\alpha(G) \leq p-k.$
\end{cor}

\begin{rem}
Setting $p \geq c \cdot k$ for a constant $c > 1$ we obtain $n_0 < p + \frac{c}{c-1}(k+1)$ in Theorem \ref{t4}. For the complexity we obtain
$O(1.2738^{2k-1+\lceil\frac{k+1}{c-1}\rceil} + kn).$
\end{rem}

%Next we show that $\frac{p(p+1)}{p-k} \leq 2pk.$ If $k=1,$ then $\frac{p(p+1)}{p-1} \leq 2p,$ since $p \geq 3.$ Hence we may assume that $k \geq 2.$
%Now $\frac{p(p+1)}{p-k} \leq 2pk \leftrightarrow 2k^2+1 \leq (2k-1)p.$ Since $2k^2+1 \leq (2k-1)(k+1)$ for $k \geq 2$ and $k \leq p-1,$ we obtain $\frac{p(p+1)}{p-k} \leq 2pk.$

%Using the vertex cover fpt-algorithm of Chen et al. \cite{CKX} we can decide in time $O(1,2738^{3k} + 3kn)$ whether $\alpha(G_0) \leq p-k.$

\section{Improving the upper bound}

Before presenting an algorithm including our main results we will show an improvement of the upper bound using neighbourhood unions of pairs
of nonadjacent vertices. This bound is easy computable and can eventually speed up the running time of the algorithm.

\medskip

We start with the following upper bound for the independence number, which may be attributed as folklore.
For a graph $G$ of order $n$ let $d_1 \leq d_2 \leq \ldots \leq d_n$ denote its degree sequence. Define
$p_1(G):=max\{i | d_i \leq n-i\}.$ Note that $p_1(G)$ is well defined, since $d_1 \leq n-1.$
Then the following proposition holds.

\begin{pro}\label{p1}
Let $G$ be a graph of order $n$ and degree sequence $d_1 \leq d_2 \leq \ldots \leq d_n.$ Then
$$\alpha(G) \leq p_1(G).$$
\end{pro}

A short proof can be deduced from the following Theorem of Welsh and Powell.

\begin{thm}(Welsh and Powell)\cite{WP}\label{t5}
Let $G$ be a graph of order $n$ and degree sequence $d_1 \geq d_2 \geq \ldots \geq d_n.$ Then
$$\chi(G) \leq \max\{\min_{1 \leq i \leq n}\{i,d_i +1\}.$$
\end{thm}

It is well known that $\alpha(G) \leq cc(G) = \chi(\overline{G})$ for every graph $G,$ where $cc(G)$ denotes the clique covering number of a graph $G.$
Now applying Theorem \ref{t5} on $\overline{G}$ we obtain $\alpha(G) \leq \max\{i | d_i \leq n-i\}.$

\medskip

Proposition \ref{p1} can be extended as follows.

\begin{pro}\label{p2}
Let $G$ be a graph of order $n$ and degree sequence $d_1 \leq d_2 \leq \ldots \leq d_n.$ Then
$$\alpha(G) \leq p_1(G) \leq p(G).$$
\end{pro}

\proof
Suppose $p_1 = p_1(G) > p(G) = p.$ Then $\sum_{i=1}^{p_1}(n-1-d_i) \geq \sum_{i=1}^{p_1}(n-1-d_{p_1}) \geq
p_1(p_1-1) \geq (p+1)p > \overline{m}(G),$ a contradiction.
\qed

\medskip

For a vertex $u \in V(G)$ and all vertices $v \in \overline{N}(u),$ let $n_2(u) \leq n_3(u) \leq \ldots \leq n_t(u)$
denote the sequence of the cardinalities of the neighbourhood unions $|N(u) \cup N(v)|$ with $t = n - d(u).$ Define
$p_2(G):= max \{k \ | \ G \ \mbox{has at least} \ k \ \mbox{vertices} \ v \ \mbox{with} \ n_k(v) \leq n-k\}.$

\begin{pro}\label{p3}
$\alpha(G) \leq p_2(G).$
\end{pro}

\proof
Let $I$ be a maximum independent set of vertices with $|I|=k.$ Let $I = \{v_1, \ldots, v_k\}.$ Then $|N(v_i) \cup N(v_j)| \leq n-k$ for all pairs of vertices $v_i, v_j$ with $1 \leq i < j \leq k.$ Hence $n_k(v_i) \leq n-k$ for $1 \leq i \leq k$ implying
$p_2(G) \geq k = \alpha(G).$
\qed

\begin{thm}
Let $G$ be a graph of order $n$ and degree sequence $d_1 \leq d_2 \leq \ldots \leq d_n.$ Then
$$\alpha(G) \leq p_2(G) \leq p_1(G) \leq p(G).$$
\end{thm}

\proof
Suppose $p_2(G) > p_1(G).$ Then $d_{p_2} > n - p_2$ by the definition of $p_1(G).$ Let $v_1, \ldots, v_{p_2}$ be $p_2(G)$ vertices with $nc_{p_2}(v_i) \leq n - p_2.$ Then $d(v_i) \leq n - p_2$
for $1 \leq i \leq p_2$ implying $d_{p_2} \leq n - p_2,$ a contradiction. This shows that $p_2(G) \leq p_1(G).$
Now the result follows with Proposition \ref{p2} and Proposition \ref{p3}. \qed

\begin{rem}
For two integers $n, p$ with $n > p \geq 2$ let $H_{n,p} = K_{n-p} + pK_1.$ Then $\alpha(H_{n,p}) = p_2(H_{n,p}) = p_1(H_{n,p}) = p(H_{n,p}).$
This shows that there are infintely many graphs $G$ with $\alpha(G) = p_2(G) = p_1(G) = p(G).$
\end{rem}

\medskip

\noindent{\bf ALGORITHM}
\begin{enumerate}
\item For a graph $G$ with $n$ vertices and $m$ edges compute $p = p(G)$ and $p_1(G).$
\item If $p_1(G) \leq p-k,$ then {\bf Output:} YES
\item Else compute $p_2(G).$
\item If $p_2(G) \leq p-k,$ then {\bf Output:} YES
\item Else compute $G_{p,k}$ and decide whether $\alpha(G_{p,k}) \leq p-k$ using the vertex cover fpt-algorithm.
\end{enumerate}

\section{Extremal graphs}

In this section we will describe the structure of extremal graphs $G_{p,k}$.

For an independent set $I$ of a graph $G$ a set $S$ is called an {\it augmenting set} for $I,$ if
$(I \setminus N(S)) \cup S$ is an independent set with $|(I \setminus N(S)) \cup S| > |I|.$

Now let $I$ be a maximum independent set in $G_0,$ and let $R = V(G_0) \setminus I.$ Then for any $t$ independent vertices
$x_1, \ldots, x_t$ in $R$ there are $t$ vertices $y_1, \ldots, y_t$ in $I$ such that $x_iy_i \in E(G)$ for $1 \leq i \leq t.$

\medskip

Suppose, $\alpha(G) \leq p - (k-1)$ and we obtain the answer "YES" to the test $\alpha(G_{p,k}) \geq p-k+1.$ Then $\alpha(G) =
\alpha(G_{p,k}) = p-k+1.$ Let $I$ be a maximum independent set in $G_{p,k},$ and set $V(G_{p,k}) = I \cup R$ with $r = |R|.$
By the construction of $G_{p,k}$ we have $d_{\overline{G}}(v) \geq p-k$ for each vertex $v \in V(G_{p,k}).$ Then

$${p+1 \choose 2} > m(\overline{G}_{p,k}) \geq |E(I)|+|E(I,R)|+|E(R)|.$$

\noindent Let $E^* = E(I,R) \cup E(R).$ Then $|E^*| = |E(I,R)|+|E(R)| \leq {p+1 \choose 2} - 1 - {p-k+1 \choose 2}.$
Hence $|E^*| \leq p-1$ for $k=1,$ $|E^*| \leq 2p-2$ for $k=2,$ and $|E^*| \leq 3p-4$ for $k=3.$  
% =
%\[ \left\{
%\begin{array}{ll}
%p-1, & if k=1, \\
%2p-2, & if k=2, \\
%3p-4, & if k=3.
%\end{array}
%\]
%$$

On the other hand, using the handshaking lemma,

$$|E^*| \geq max\{r(p-k) - {r \choose 2}, \frac{r(p-k)}{2}\}.$$

\noindent This inequality has the following properties.

\medskip

\noindent (P1) \quad $r(p-k) - {r \choose 2}$ is increasing for $1 \leq r \leq p-k.$

\medskip

\noindent (P2) \quad $\frac{r(p-k)}{2} \geq (p-k+1)\frac{p-k}{2}$ for $r \geq p-k+1.$

\medskip

Now the extremal graphs can be generated, if we require $p \geq p(k)$ for a given integer $k \geq 0.$

\begin{enumerate}
\item
$\alpha(G) = p, k=1, p \geq 3$ \\

We have $\frac{r(p-1)}{2} \leq |E^*| \leq p-1,$ hence $r \leq 2.$ By (P1) we
have $2(p-1)-1 = 2p-3 > p-1$ for $r=2$ and $p \geq 3.$ Hence we may assume $0 \leq r \leq 1$ and we distinguish these two cases. Let $I$ be a maximum independent set of size $|I|=p.$ Then there are no augmenting sets for $I.$ This leads to the following graph structures.

\begin{enumerate}
\item
$n_0 = p, r=0$\\
Then $G_{p,1} \cong pK_1.$
\item
$n_0 = p+1, r=1$\\
Then $K_2 \cup (p-1)K_1 \subset G_{p,1} \subset pK_1 + K_1.$
\end{enumerate}

\item
$\alpha(G) = p-1, k=2, p \geq 8$ \\

We have $\frac{r(p-2)}{2} \leq |E^*| \leq 2p-2.$
If $r \geq p-k+1=p-1,$ then by (P2) we obtain $|E^*| \geq (p-1)\frac{5}{2} > 2p-2$ for $p \geq 7,$ a contradiction.
Hence we have $r \leq p-2.$
By (P1) we
have $3(p-2)-3 = 3p-9 > 2p-2$ for $r=3$ and $p \geq 8.$ Hence we may assume $0 \leq r \leq 2$ and we distinguish these three cases.
Let $I$ be a maximum independent set of size $|I|=p-1.$ Then there are no augmenting sets for $I.$ This leads to the following graph structures.

%Then $p \leq n_0 = |V(G_0(p))| \leq p+1.$ ????
\begin{enumerate}
\item
$n_0 = p-1, r=0$\\
Then $G_{p,2} \cong (p-1)K_1.$
\item
$n_0 = p, r=1$\\
Then $K_2 \cup (p-2)K_1 \subset G_{p,2} \subset (p-1)K_1 + K_1.$
\item
$n_0 = p+1, r=2$\\
Then $K_3 \cup (p-2)K_1 \subset G_{p,2} \subset (p-1)K_1 + K_2$\\ or \quad
$2K_2 \cup (p-3)K_1 \subset G_{p,2} \subset (p-1)K_1 + K_2$
\end{enumerate}

\item
$\alpha(G) = p-2, k=3, p \geq 15$ \\

We have $\frac{r(p-2)}{2} \leq |E^*| \leq 3p-4.$
If $r \geq p-k+1=p-2,$ then by (P2) we obtain $|E^*| \geq (p-2)\frac{p-3}{2} > 3p-4$ for $p \geq 10,$ a contradiction.
Hence we have $r \leq p-3.$
By (P1) we
have $4(p-3)-6 = 4p-18 > 3p-4$ for $r=4$ and $p \geq 15.$ Hence we may assume $0 \leq r \leq 3$ and we distinguish these four cases.
Let $I$ be a maximum independent set of size $|I|=p-2.$ Then there are no augmenting sets for $I.$ This leads to the following graph structures.

%Then $p \leq n_0 = |V(G_0(p))| \leq p+1.$ ????
\begin{enumerate}
\item
$n_0 = p-2, r=0$\\
Then $G_{p,3} \cong (p-2)K_1.$
\item
$n_0 = p-1, r=1$\\
Then $K_2 \cup (p-3)K_1 \subset G_{p,3} \subset (p-2)K_1 + K_1.$
\item
$n_0 = p, r=2$\\
Then $K_3 \cup (p-3)K_1 \subset G_{p,3} \subset (p-2)K_1 + K_2$\\ or \quad
$2K_2 \cup (p-4)K_1 \subset G_{p,3} \subset (p-2)K_1 + K_2$
\item
$n_0 = p+1, r=3$\\
Then $K_4 \cup (p-3)K_1 \subset G_{p,3} \subset (p-2)K_1 + K_3$\\
or \quad
$K_3 \cup K_2 \cup (p-4)K_1 \subset G_{p,3} \subset (p-2)K_1 + K_3$\\
or \quad
$3K_2 \cup (p-5)K_1 \subset G_{p,3} \subset (p-2)K_1 + K_3$
\end{enumerate}

\end{enumerate}

%\bigskip

%\newpage

\end{document}